\documentclass[reqno]{amsart}
\usepackage{amscd}

\def\beq{\begin{equation}}
\def\eeq{\end{equation}}
\def\ba{\begin{array}}
\def\ea{\end{array}}
\def\S{\mathbb S}
\def\R{\mathbb R}

\def\la{\langle}
\def\ra{\rangle}

\def \ds{\displaystyle}
\def \vs{\vspace*{0.1cm}}

\def\slashii#1{\setbox0=\hbox{$#1$}             
   \dimen0=\wd0                                 
   \setbox1=\hbox{\sl/} \dimen1=\wd1            
   \ifdim\dimen0>\dimen1                        
      \rlap{\hbox to \dimen0{\hfil\sl/\hfil}}   
      #1                                        
   \else                                        
      \rlap{\hbox to \dimen1{\hfil$#1$\hfil}}   
      \hbox{\sl/}                               
   \fi}                                         %
\def\slashiii#1{\setbox0=\hbox{$#1$}#1\hskip-\wd0\hbox to\wd0{\hss\sl/\/\hss}}
\newcommand{\C}{{\mathbf C}}

\newtheorem{thm}{Theorem}[section]
\newtheorem{lm}[thm]{Lemma}
\newtheorem{prop}[thm]{Proposition}

\theoremstyle{definition}
\newtheorem{rem}[thm]{Remark}

\theoremstyle{remark}

\begin{document}
\pagestyle{plain}
\today

\title{Super-Liouville Equations on Closed Riemann Surfaces}
\author{J\"urgen Jost}
\address{Max Planck Institute for Mathematics in the Sciences, Inselstr. 22, D-04013 Leipzig}
\email{jjost@mis.mpg.de}
\author{Guofang Wang}
\address{Max Planck Institute for Mathematics in the Sciences, Inselstr. 22, D-04013 Leipzig}
\email{gwang@mis.mpg.de}
\author{Chunqin Zhou}
\address{Department of Mathematics, Shanghai Jiaotong University, Shanghai, 200240, China }
\email{cqzhou@sjtu.edu.cn}
\thanks{The third named author
supported partially by NSFC of China (No. 10301020)}

\begin{abstract}
Motivated by the supersymmetric extension of Liouville theory in
the recent physics literature, we couple the standard Liouville
functional with  a spinor field term. The resulting functional is
conformally invariant. We study  geometric and analytic aspects of
the resulting Euler-Lagrange equations, culminating in a blow up
analysis.
\end{abstract}
\maketitle

\section{Introduction}
The classical Liouville functional  for a real-valued function $u$ on $M$ is
\begin{equation}
\label{functional1} E\left( u \right) =\int_{M}\{\frac 12 \left|
\nabla u\right| ^2+K_gu -e^{2u}\}dv,
\end{equation}
where  $K_g$ is the Gaussian curvature of $M$.
 The Euler-Lagrange equation for $E(u)$ is the Liouville equation
\begin{equation}
\label{1.1}
-\Delta u = \ds\vs 2e^{2u}  -K_g
\end{equation}
where $\Delta $ is the Laplacian with respect to $g$. Liouville
\cite{L} studied this equation in the plane, that is, for $K_g=0$.
The Liouville equation arises in many contexts of complex analysis
and differential geometry of Riemann surfaces, in particular in
the prescribing curvature problem. The interplay between the
geometric and analytic aspects makes the Liouville equation
mathematically rich. It also occurs naturally in string theory as
discovered by Polyakov \cite{Po2}, from the gauge anomaly in
quantizing the string action. There then also is a natural
supersymmetric version of the Liouville functional and equation,
coupling the bosonic scalar field to a fermionic spinor field. It
turns out, however, that we also obtain a very interesting
mathematical structure if we consider ordinary instead of
fermionic (Grassmann valued) spinor fields. In particular, the
fundamental conformal invariance of the Liouville action can be
preserved under the coupling. This makes the resulting functional
geometrically very natural and, so it seems to us, a worthy and
interesting object of mathematical analysis.

Therefore, in this paper, we consider
the following functional for a real-valued function $u$ and a spinor $\psi$
\begin{equation}
\label{functional2} E\left( u,\psi \right) =\int_{M}\{\frac 12
\left| \nabla u\right| ^2+K_gu+\left\langle (\slashiii{D}+e^u)\psi
,\psi \right\rangle -e^{2u}\}dv.
\end{equation}
The Euler-Lagrange system for $E(u,\psi )$ is
\begin{equation}
\left\{
\begin{array}{rcl}
-\Delta u &=& \ds\vs 2e^{2u}-e^u\left\langle \psi ,\psi
\right\rangle  -K_g\qquad
\\
\slashiii{D}\psi &=&\ds  -e^u\psi
\end{array}
\text{in }M,\right.   \label{1}
\end{equation}
This system couples  the Liouville equation
and the Dirac equation in a rather natural way.
We  call (\ref{1}) the {\it super-Liouville equations}.
 When $\psi $ vanishes, we obtain the original Liouville equation.  In
other words, here we are considering a system generalizing the
prescribing curvature equation. The important point is that this
generalization preserves a fundamental property of the energy
functional on Riemann surfaces, namely its conformal invariance.

In this paper we aim to provide an analytic foundation for system
(\ref{1}). We start with basic points like the regularity of weak
solutions. An analytic foundation was established for the
Liouville equation (\ref{1.1}) in \cite{BM},  \cite{LS} and for a
Toda system  in \cite{JW}, \cite{JW2} and \cite{JLW}. In those
references, it was established that the key analytical points are
that singularities in solutions $u_n$ of the equations on closed
surfaces, or, more generally with bounded energy $\int e^{2u_n}$,
can form only at isolated points $x$ where the limit $u_n(x)$
tends to infinity. Away from those singularities, $u_n$ remains
either uniformly bounded or converges to $-\infty$ which, in fact,
is a regular situation for the field $\phi$ with $u=\log \phi$. At
those isolated singularities, rescaling produces an entire
solution of the Liouville equation of finite  energy
$\int_{\R^2}e^{2u}$ in the  plane which then can be compactified
to a solution on the 2-sphere. Therefore, the asymptotic behavior
of such entire solutions is also an important point. In this
paper, we therefore perform such an analysis for the super
Liouville equations. As in the classical case, this provides a
complete analytical picture, and other regularity results follow
in a standard manner that is known to the experts and therefore
need not be repeated here.

Assume that $(u_n,\psi _n)$ is a sequence of solutions of (\ref{1}) with
\[
\int_{M}e^{2u_n}dv<\varepsilon _0,\text{ and }\int_{M}\left| \psi
_n\right| ^4dv<C
\]
for some positive constants $\varepsilon _0$ and $C$. If
$\varepsilon _0$ is sufficiently small (in fact $\varepsilon
_0<\pi$ suffices), then we can show that $(u_n,\psi _n)$ admits a
subsequence, which we still denote by $(u_n, \psi_n)$, converging
to a smooth solution $(u,\psi)$ of (\ref{1}). Note that
$\int_{M}e^{2u_n}dv$ and $ \int_{M}\left| \psi _n\right| ^4dv$ are
conformally invariant, see Section 3.

When $\varepsilon_0$ is big, then the so-called ``blow-up" phenomenon may occur.
Let ($u_n,\psi _n)$ be a sequence of solutions of
(\ref{1}) and satisfying
\[
\int_{M}e^{2u_n}dv<C,\text{ and }\int_{M}\left| \psi _n\right|
^4dv<C.
\]
Define

\[\ba{rcl}
\ds \Sigma _1 &=&\ds\vs \left\{ x\in M,\text{ there is a sequence
}y_n\rightarrow x\text{ such that }u_n(y_n)\rightarrow +\infty
\right\}\\
\ds \Sigma _2&=& \ds\left\{ x\in M,\text{ there is a sequence
}y_n\rightarrow x\text{ such that }\left| \psi _n(y_n)\right|
\rightarrow +\infty \right\}.\ea
\]
Then, one can show that $\Sigma_2\subset \Sigma_1$ and
 ($u_n,\psi _n)$ admits a subsequence, still denoted by ($%
u_n,\psi _n),$ satisfying one of the following cases:

\begin{enumerate}
\item[i)]  $u_n$ is bounded in $L^\infty (M).$

\item[ii)]  $u_n$ $\rightarrow -\infty $ uniformly on $M$.

\item[iii)]  $\Sigma _1$ is finite, nonempty and either
\[
u_n\text{ is bounded in }L_{loc}^\infty (M\backslash \Sigma _1)
\]
or
\[
u_n\rightarrow -\infty \text{ uniformly on compact subsets of
}M\backslash \Sigma _1.
\]
\end{enumerate}

Furthermore, we rule out the first case in iii) if $\Sigma
_1\backslash \Sigma_2 \neq \emptyset $. Then the only case is
$u_n\rightarrow -\infty $ uniformly on compact subsets of
$M\backslash \Sigma _1.$

Finally, we consider entire solutions of the super Liouville equations on
 $\R^2$ with finite energy $\int_{\R^2}e^{2u}+|\psi|^4$, which
 can  be viewed as ``bubbles"
or obstructions to the compactness of equation (\ref{1}).
 We analyze the asymptotic behavior of  such solutions
and obtain
\begin{equation*}
u(x)=-\frac {\alpha}{2\pi} \ln{|x|}+C+O(|x|^{-1}) \qquad
\text{for}\quad |x| \quad \text{near}\quad \infty,
\end{equation*}

\begin{equation*}
\psi (x)=-\frac {1}{2\pi}\frac{x}{|x|^2}\cdot
\xi_0+o(|x|^{-1})\qquad \text{for}\quad |x| \quad
\text{near}\quad \infty,
\end{equation*}
where $\cdot$ is the Clifford multiplication, $C\in \R$ is some
constant, $\alpha =\int_{\R^2}2e^{2u}- e^u|\psi|^2dx$, and
$\xi_0=\int_{\R^2}e^u\psi dx$ is a constant spinor.

Furthermore,  by using the associated holomorphic quantity
$T(z)=(\partial_z u)^2-\partial^{2}_{z}u+\frac 14\langle
\psi,dz\cdot\partial_{\bar z}\psi\rangle +\frac 14\langle
d\bar{z}\cdot\partial_z\psi,\psi\rangle$, we show $\alpha =4\pi$.
For the definition of $T$, see Section 3.
 Then we show  such an entire solution
can be extended to a smooth solution on $\S^2$, i.e. the  global
singularity (the singularity at infinity)
is removable.

\section{Spinors}
For presenting our equations, we need to recall some background about spin structures and spinors.
Let $(M,g)$ be a  closed Riemann surface and $P_{SO(2)}\to M$ its
oriented orthonormal frame bundle. A $Spin$-structure is a lift of
the structure group SO(2) to $Spin(2)$, {\it i.e.}, there exists a
principal $Spin$-bundle $P_{Spin(2)}\to M$ such that there is a
bundle map
\[ \begin{array}{ccc}
 P_{Spin(2)} & \longrightarrow & P_{SO(2)} \cr
&& \cr \downarrow & & \downarrow \cr && \cr M & \longrightarrow &
M.\cr
\end{array}\]

Let $\Sigma ^+M:=P_{Spin(2)}\times_{\rho}{\mathbb C}$ be a complex
line bundle over $M$ associated to $ P_{Spin(2)}$ and to the
standard representation $\rho:\S^1\to U(1)$. This is the bundle of
positive half-spinors. Its complex conjugate $\Sigma^-M:=
\overline {\Sigma^+M}$ is called the bundle of negative
half-spinors. The spinor bundle is $\Sigma M:=\Sigma ^+M\oplus
\Sigma ^-M.$
There exists a Clifford multiplication
\begin{eqnarray*}
TX\times_{\C} \Sigma^+M &\to &  \Sigma^-M \cr TX\times_{\C}
\Sigma^-M &\to &  \Sigma^+M \end{eqnarray*} denoted by $v\otimes
\psi \to v \cdot \psi$, which satisfies the Clifford relations
\[ v\cdot w \cdot\psi + w\cdot v\cdot \psi =-2g(v,w) \psi,\]
for all $v,w \in TM$ and $\psi \in \Gamma(\Sigma M)$.

On the spinor bundle $\Sigma M$, the metric $g$ induces a natural Hermitian
metric $ \la \cdot,\cdot \ra $. Let $\nabla$ be
the Levi-Civita connection on $M$ with respect to $g$. Likewise, $\nabla$ induces a
connection (also denoted by $\nabla$) on  $\Sigma M$ compatible
with the Hermitian metric.

The Dirac operator $\slashiii{D}$ is defined by $\slashiii{D}\psi
:=\sum_{\alpha =1}^2e_{_\alpha }\cdot \nabla _{e_\alpha }\psi ,$
where $\left\{ e_1,e_2\right\} $ is an orthonormal basis on $TM$.
(For more details about the spin bundle and the Dirac operator, we
refer to \cite{LM} or \cite{Jo}.)

\section{Properties of Super Liouville Equations}
We start by giving some examples of solutions of the super Liouville equations (\ref{1}).
When $M=\S^2$, the standard sphere with Gaussian curvature $K=1$,
it is obvious that   solutions $u$ of (\ref{1.1}),
\begin{equation}
-\Delta u+1-2e^{2u}=0 \text{ on } \S^2.  \label{2}
\end{equation}
  yield solutions of the form $(
u,0)$ of (\ref{1}),
\begin{equation}
\left\{
\begin{array}{rcl}
-\Delta u &=& \ds\vs 2e^{2u}-e^u\left\langle \psi ,\psi
\right\rangle  -1\qquad
\\
\slashiii{D}\psi &=&\ds  -e^u\psi
\end{array}
\right.
\end{equation}
In fact, all solutions of (\ref{2}) are of the form $u=\frac 12
\log \frac 12+\frac 12\log \det \left| d\varphi \right| ,$ where
$\varphi $ is a conformal map of $\S^2$. This can be understood in
terms of the complex geometry behind the Liouville equation, but
we do not go into this aspect here.

There  exists another type of solution of (\ref{1}). Let us recall
that a {\it Killing spinor} is a spinor $\psi$  satisfying
\[
\nabla _X\psi = \lambda X\cdot \psi , \quad \hbox{
for any vector field } X
\]
for some constant $\lambda$. On the standard sphere, there are Killing spinors
with the Killing constant $\lambda= \frac 12$, see for instance \cite{baum}.
 Such a
Killing spinor
is an eigenspinor, i.e.
\[
\slashiii{D}\psi =- \psi,
\]
with constant $|\psi|^2$. Choosing a Killing spinor $\psi$ with
$|\psi|^2=1$, $(0,\psi)$ is a solution
of (\ref{1}). If we identify $\S^2\backslash\{north pole\}$ by the
stereographic projection with the Euclidean plane $\R^2$ with the
metric
\[\frac 4{(|1+|x|^2)^2}|dx|^2,\]
then any Killing spinor has the form
\[\frac{v+x\cdot v}{\sqrt{1+|x|^2}},\]
up to a translation or a dilation. See \cite{baum}.\\\\

Now we come to an important property of the
functional $E$.

\begin{prop}\label{prop-c}
The functional $E(u,\psi )$ is conformally invariant. Namely, for
any conformal diffeomorphism $\varphi :M\rightarrow M,$ set
\beq\label{n3.1}\ba{rcl}
\ds\vs \widetilde{u} &=&\ds u\circ \varphi -\ln \lambda  \\
\ds \widetilde{\psi } &=&\ds \lambda ^{-\frac 12}\psi \circ \varphi \ea
\eeq
where $\lambda $ is
the conformal factor of the conformal map $\varphi,$ i.e.,
$\varphi^*(g)=\lambda ^2 g$. Then $E(u,\psi )=E(\widetilde{u},\widetilde{\psi })$.
In particular, if $(u,\psi)$ is a solution of
(\ref{1}), so is $(\tilde u,\tilde \psi)$.
\end{prop}
\begin{proof}
It is well-known that $\int_{M}\frac 12\left| \nabla u\right|
^2+K_gu$ is conformally invariant, see e.g. \cite{H}. Since the terms
\[
\int_{M}e^{2u}dv,\qquad \int_{M}e^u\left| \psi \right| ^2dv,
\]
are invariant under a conformal transformation, it is sufficient
to show the
conformality of $\int_{M}\left\langle \slashiii{D}\psi ,\psi \right\rangle dv.$ Let $%
\widetilde{g}=\varphi ^{*}g$, where $g$ is the metric on $M.$ Let
$\widetilde{\slashiii{D}}$ be the Dirac operator with respect to
the new
metric $\widetilde{g}.$ By the conformality of $\varphi ,$ we have $%
\widetilde{g}=\lambda ^2g$ for a positive function $\lambda $ on
$M$. We identify the new and old spin bundles as in \cite{H}.
Since the relation between the two Dirac operators $\slashiii{D}$ and
$\widetilde{\slashiii{D}}$ is
\[
\widetilde{\slashiii{D}}\widetilde{\psi }=\lambda ^{-\frac 32}\slashiii{D}(\lambda ^{\frac 12}%
\widetilde{\psi })=\lambda ^{-\frac 32}\slashiii{D}\psi ,
\]
we can show by a direct computation that
\[
\int_{M}\left\langle \slashiii{D}\psi ,\psi \right\rangle
dvol(g)=\int_{M}\langle \widetilde{\slashiii{D}}\widetilde{\psi },\widetilde{%
\psi }\rangle dvol(\widetilde{g}).
\]
The proof of the proposition is complete.%
\end{proof}

As before, we identify $\S^2\backslash\{north pole\}$ by
stereographic projection with the Euclidean plane $\R^2$ with the
metric
\[\frac 4{(|1+|x|^2)^2}|dx|^2.\]

By Proposition \ref{prop-c}  from any solution of equation
(\ref{1}) on $\S^2$ one can obtain a solution of

\begin{equation}\label{eqR2}
\left\{
\begin{array}{rcl}
-\Delta u &=& \ds\vs 2e^{2u}-e^u\left\langle \psi ,\psi
\right\rangle
\\
\slashiii{D}\psi &=&\ds  -e^u\psi\qquad
\end{array}
\right.
\text{in }\R^2,
\end{equation}
where $\Delta$ and $ \slashiii{D}$ are operators with respect to the
standard metric on $\R^2$.

Equation (\ref{eqR2}) is very interesting, since its solutions are
 obstructions for the compactness of equation (\ref{1}), namely
they are the so-called
``bubbles" in the geometric analysis.

Let us note that on a surface the (usual) Dirac operator
$\slashiii{D} $ can be seen as the (doubled) Cauchy-Riemann operator.
Consider $\R^2$ with the Euclidean metric $dx^2_1+dx^2_2$. Let
$e_1=\frac{\partial}{\partial x_1}$ and
$e_2=\frac{\partial}{\partial x_2}$ be the standard orthonormal
frame. A spinor field is simply a map $\Psi:\R^2\to
\Delta_2=\C^2$, and $e_1$ and $e_2$ acting on  spinor fields can
be identified by multiplication with matrices
\[e_1=\left(\begin{matrix}0& 1\\ -1&0 \end{matrix}\right),
\quad e_2=\left(\begin{matrix}0& i\\ i&0 \end{matrix}\right).\] If
$\Psi:=\left(\begin{matrix} \ds f
\\ \ds g\end{matrix}\right)
:\R^2\to \C^2$ is a spinor field, then the Dirac operator is
\[\slashiii{D}\Psi=\ds \left(\begin{matrix}0& 1\\ -1&0
\end{matrix}\right) \left(\begin{matrix} \ds \frac{\partial
f}{\partial x_1}\\ \ds \frac{\partial g}{\partial x_1}
\end{matrix}\right)+
\left(\begin{matrix}0& i\\ i&0 \end{matrix}\right)
\left(\begin{matrix} \ds \frac{\partial f}{\partial x_2} \\
\ds\frac{\partial g}{\partial x_2}
\end{matrix}\right)=
2\left(\begin{matrix} \ds \frac{\partial g}{\partial \bar z}
\\ -\ds\frac{\partial f}{\partial z}\end{matrix}\right),\]
where
\[\frac{\partial}{\partial z}=\frac 12 \left(\frac{\partial }{\partial x_1}
- i\frac{\partial }{\partial x_2}\right), \quad
\frac{\partial}{\partial \bar z}=\frac 12 \left(\frac{\partial
}{\partial x_1} + i\frac{\partial }{\partial x_2}\right).\]
Therefore, the elliptic estimates developed for (anti-)
holomorphic functions can be used to study the Dirac equation.

\begin{prop}
\label{prop2}Let $M=\S^2$ and $\psi$
 a Killing spinor with $|\psi|=1$. Then
\[ (0, \psi)\]
 is a solution of (\ref{1}),
\end{prop}

\begin{proof}
This is obvious, and we have observed it above. In order to
understand the conformal invariance of the super Liouville
equation better, it is instructive to carry out the proof  on
$\R^2$. From the above discussion and Prop.\ref{prop-c}, it is
sufficient to show that
\begin{equation}\left\{\ba {rcl}\label{solu}
u &=& \ds\vs -\log (1+|x|^2)+\log 2, \\
\psi &=& \ds(\frac 2{|1+|x|^2})^{\frac 12}\frac{v+x\cdot
v}{\sqrt{1+|x|^2}}
\ea\right.
\end{equation}
with $v\in \{v\in {\mathbb C}^2\,|\, |v|=1\}$  is a solution of equation (\ref{eqR2}).

We write $x\cdot v=x_1e_1\cdot v+x_2e_2\cdot v$. Recall the
Clifford multiplication
\[
e_i\cdot e_j\cdot \psi+e_j\cdot e_i\cdot \psi=-2\delta
_{ij}\psi,\text{ for }1\leq i,j\leq 2
\] and
\[
\left\langle \psi ,\varphi \right\rangle =\left\langle e_i\cdot
\psi ,e_i\cdot \varphi \right\rangle
\]
for any spinor fields  $\psi ,\varphi $. It is clear that
$$\la x\cdot v,
x\cdot v \ra=|x|^2, \qquad \text{and}\qquad \la v, x\cdot v\ra+\la
x\cdot v, v\ra=0.$$

\noindent Then by a direct computation, we have
\begin{eqnarray*}
\la \psi, \psi \ra  &=&  \frac 2{(1+|x|^2)^2}\la v+x\cdot v,
v+x\cdot
v\ra\\
&=&  \frac 2{(1+|x|^2)^2}(\la v, v\ra + \la x\cdot v, x\cdot v \ra + \la v, x\cdot v\ra+\la x\cdot v, v\ra )\\
&=& \frac 2{1+|x|^2}.
\end{eqnarray*}
 Thus we can easily check that
$(u,\psi)$ satisfies the first equation.

Next we calculate that
$$
\partial_{x_1}\psi=\frac {-2\sqrt{2}x_1}{(1+|x|^2)^2}(v+x\cdot
v)+\frac {\sqrt{2}}{1+|x|^2}e_1\cdot v,
$$
and
$$
\partial_{x_2}\psi=\frac {-2\sqrt{2}x_2}{(1+|x|^2)^2}(v+x\cdot
v)+\frac {\sqrt{2}}{1+|x|^2}e_2\cdot v.
$$

\noindent Then we have
\begin{eqnarray*}
\slashiii{D}\psi &=& e_1\cdot \partial_{x_1}\psi+e_2\cdot
\partial_{x_2}\psi\\
&=& -\frac {2\sqrt{2}x_1}{(1+|x|^2)^2}(e_1\cdot v-x_1
v+x_2e_1\cdot e_2\cdot v)\\
& &-\frac {2\sqrt{2}x_2}{(1+|x|^2)^2}(e_2\cdot v-x_2 v+x_1e_2\cdot
e_1\cdot v)-\frac
{2\sqrt{2}}{1+|x|^2} v\\
&=& -\frac {2\sqrt{2}}{(1+|x|^2)^2}(v+x\cdot v) \\
&=& -e^u\psi.
\end{eqnarray*}
This implies that $(u,\psi)$ satisfy the second equation.
\end{proof}

\
\

By  conformal transformations, we know that
\[
(\log\frac {\sqrt{2}}{1+|x-x_0|^2},0) \hbox{  and } (\log\frac
2{1+|x-x_0|^2}, \sqrt {2}\frac {v+(x-x_0)\cdot v}{1+|x-x_0|^2})\]
are solutions of (\ref{eqR2}). It is clear that all such solutions
of (\ref{eqR2}) obtained from solutions of (\ref{1}) on $\S^2$
satisfy
\begin{equation}\label{bound}
I(u,\psi):=\int_{\R^2}\{|\nabla u|^2+|\psi|^4\}<C.
\end{equation}
In the last section, we will show that all solutions of
(\ref{eqR2}) with bounded energy $I$ are obtained from solutions
of (\ref{1}) on $\S^2$.

\begin{prop}\label{pro2}
Let $(u,\psi )$ be a smooth solution of (\ref{1}) and $z=x+iy$  a
local isothermal parameter with $g=ds^2=\rho \left| dz\right|^2 .$
Then the quadratic differential
\[
T(z)dz^2=\{(\partial_z u)^2-\partial^{2}_{z}u+\frac 14\langle
\psi,dz\cdot\partial_{\bar z}\psi\rangle +\frac 14\langle
d\bar{z}\cdot\partial_z\psi,\psi\rangle  \}dz^2
\]
is holomorphic when $M$ is a constant curvature surface. Here
$dz=dx+idy$ and $d\bar z =dx-idy$.
\end{prop}

\begin{proof}
We prove this lemma by a direct computation. Let  \{$e_1,e_2\}$ be
a local orthonormal basis on $M$. It follows from the Clifford
multiplication  that
$$
\langle e_\alpha\cdot\psi,\psi\rangle=\langle e_\alpha\cdot
e_\alpha\cdot\psi, e_\alpha\cdot\psi\rangle=-\langle \psi,
e_\alpha\cdot\psi\rangle.
$$
Therefore we obtain the real part of $\langle
e_\alpha\cdot\psi,\psi\rangle$  vanishes, i.e. \beq\label{3.1}
\text{ Re }\langle e_\alpha\cdot\psi,\psi\rangle=0. \eeq
Furthermore we have
\[
\left\langle \psi ,e_{\alpha }\cdot \nabla _{e_{\beta }}\psi
\right\rangle =-\left\langle e_{\alpha }\cdot \psi ,\nabla
_{e_{\beta }}\psi \right\rangle =-\left\langle e_{\beta }\cdot
e_{\alpha }\cdot \psi ,e_{\beta }\cdot \nabla _{e_{\beta }}\psi
\right\rangle
\]
and
\begin{eqnarray*}
\left\langle \psi ,e_{\alpha }\cdot \nabla _{e_{\beta }}\psi
\right\rangle -\left\langle \psi ,e_{\beta }\cdot \nabla
_{e_{\alpha }}\psi \right\rangle &=&-\left\langle e_{\beta }\cdot
e_{\alpha }\cdot \psi
,\slashiii{D}\psi \right\rangle \\
&=&\left\langle e_{\beta }\cdot e_{\alpha }\cdot \psi ,e^u\psi
\right\rangle,
\end{eqnarray*}
Hence from (\ref{3.1}) we have $\text{Re}\left\langle \psi
,e_{\alpha }\cdot \nabla _{e_{\beta }}\psi \right\rangle $ is
symmetric.

\noindent Set
$$
T_1(z)=(\partial_zu)^2-\partial_z^2u,
$$
and
$$
T_2(z)=\langle \psi , dz\cdot \partial_{\bar z}\psi
\rangle+\langle d\bar{z}\cdot \partial_z\psi,\psi \rangle.
$$

\noindent Then, we choose a local orthonormal basis \{$e_1,e_2\}$
on $M$ such that $\nabla _{e_{\alpha }}e_{\beta }=0$ at a
considered point.
By using the Ricci curvature formula we have
$$
\partial_{zz\overline{z}}u=\frac 14(\partial_z (\triangle u)+2K_g\partial_z
u).
$$
Now we can compute
\begin{eqnarray*}
\partial_{\bar{z}}T_1(z)&=&
2\partial_{z\overline{z}}u\partial_z u-\partial_{zz\overline{z}}u\\
&=&\frac 12 \triangle u \partial_z u-\frac 14
\partial_{z}(\triangle u)-\frac 12 K_g\partial_z u\\
&=&\frac 12 (-2e^{2u}+ e^u|\psi|^2+K_g)\partial_z u+\frac 14
\partial_z (2e^{2u}-e^u|\psi|^2-K_g)-\frac 12 K_g\partial_z u\\
&=&\frac 14e^u\left| \psi \right| ^2\partial_z u-\frac
14e^u\partial_z\left| \psi \right|^2-\frac 14
\partial_zK_g.
\end{eqnarray*}

\noindent By using the symmetry of $\text{Re}\left\langle \psi
,e_{\alpha }\cdot \nabla _{e_{\beta }}\psi \right\rangle $, we
have
\begin{eqnarray*}
\partial_{\bar{z}} T_2(z)&=&\frac
12\partial_{\bar{z}}(\langle(e_1-ie_2)\cdot(\nabla_{e_1}\psi-i\nabla_{e_2}\psi),\psi\rangle+
\langle\psi, (e_1+ie_2)\cdot(\nabla_{e_1}\psi+i\nabla_{e_2}\psi)\rangle)\\
&=&
\partial_{\bar{z}}(\text{Re}\langle\psi,e_1\cdot\nabla_{e_1}\psi\rangle
-2i\text{Re}\langle\psi,e_1\cdot\nabla_{e_2}\psi\rangle-\text{Re}\langle\psi,e_2\cdot\nabla_{e_2}\psi\rangle)\\
&=& \frac 12 (\text{Re}\langle \nabla
_{e_1}\psi,e_1\cdot\nabla_{e_1}\psi\rangle-2i\text{Re}\langle
\nabla _{e_1}\psi,e_1\cdot\nabla_{e_2}\psi\rangle-\text{Re}\langle
\nabla
_{e_1}\psi,e_2\cdot\nabla_{e_2}\psi\rangle)\\
&+&\frac 12(i\text{Re}\langle \nabla
_{e_2}\psi,e_1\cdot\nabla_{e_1}\psi\rangle+2\text{Re}\langle
\nabla
_{e_2}\psi,e_2\cdot\nabla_{e_1}\psi\rangle-i\text{Re}\langle
\nabla
_{e_2}\psi,e_2\cdot\nabla_{e_2}\psi\rangle)\\
&+&\frac 12 (\text{Re}\langle
\psi,e_1\cdot\nabla_{e_1}\nabla_{e_1}\psi\rangle-2i\text{Re}\langle
\psi,e_1\cdot\nabla_{e_1}\nabla_{e_2}\psi\rangle-\text{Re}\langle
\psi,e_2\cdot\nabla_{e_1}\nabla_{e_2}\psi\rangle)\\
&+&\frac 12 (i\text{Re}\langle
\psi,e_1\cdot\nabla_{e_2}\nabla_{e_1}\psi\rangle+2\text{Re}\langle
\psi,e_2\cdot\nabla_{e_2}\nabla_{e_1}\psi\rangle-i\text{Re}\langle
\psi,e_2\cdot\nabla_{e_2}\nabla_{e_2}\psi\rangle);
\end{eqnarray*}

\noindent It follows from (\ref{3.1}) that
$$ \text{Re}\langle \nabla
_{e_i}\psi,e_j\cdot\nabla_{e_i}\psi\rangle=0
$$
for any $i,j=1,2 $. Furthermore, by using the definition of  the
curvature operator $R^{\Sigma M}$  of the connection $\nabla $ on
the spinor bundle $\Sigma M, $ that is
\[
\nabla _{e_{\alpha} }\nabla _{e_{\beta }}\psi -\nabla _{e_{\beta}
}\nabla _{e_{\alpha }}\psi =R^{\Sigma M}(e_{\alpha },e_{\beta}
)\psi ,
\]
and a formula for this curvature operator (see for example \cite{Jo})
\[
\sum_{\alpha =1}^2e_{\alpha} \cdot R^{\Sigma M}(e_{\alpha} ,X)\psi
=\frac 12Ric(X)\cdot \psi ,\text{ \qquad for }\forall X\in \Gamma
(TM)
\]
we can obtain that
\begin{eqnarray*}
\partial_{\bar{z}} T_2(z)&=& \frac 12(-3\text{Re}\langle \nabla
_{e_1}\psi,e_2\cdot\nabla_{e_2}\psi\rangle+3i\text{Re}\langle
\nabla
_{e_2}\psi,e_1\cdot\nabla_{e_1}\psi\rangle)\\
& & +\frac 12(\text{Re}\langle
\psi,\nabla_{e_1}(\slashiii{D}\psi)\rangle-i\text{Re}\langle
\psi,\nabla_{e_2}(\slashiii{D}\psi)\rangle)\\
& & + (\text{Re}\langle \psi,e_2\cdot R^{\Sigma M
}(e_2,e_1)\psi\rangle-i\text{Re}\langle \psi,e_1\cdot R^{\Sigma M
}(e_1,e_2)\psi\rangle).
\end{eqnarray*}

\noindent By (\ref{3.1}) we have
$$
\text{Re}\langle \psi,e_2\cdot R^{\Sigma M
}(e_1,e_2)\psi\rangle=\text{Re}\langle \psi,\frac 12
Ric(e_1)\cdot\psi\rangle=0,
$$
and
$$
\text{Re}\langle \psi,e_1\cdot R^{\Sigma M
}(e_1,e_2)\psi\rangle=\text{Re}\langle \psi,\frac 12
Ric(e_2)\cdot\psi\rangle=0.
$$
We also have
\begin{eqnarray*}
\text{Re}\langle \nabla _{e_1}\psi,e_2\cdot\nabla_{e_2}\psi\rangle
\ds\vs &=& \ds \text{Re}\langle \nabla _{e_1}\psi,-e^u\psi -
e_1\cdot\nabla_{e_1}\psi\rangle\\
&=&\ds\vs \text{Re}\langle \nabla
_{e_1}\psi,-e^u\psi\rangle-\text{Re}\langle \nabla
_{e_1}\psi,e_1\cdot\nabla_{e_1}\psi\rangle\\
&=&\ds - \frac 12 e^u\nabla_{e_1}|\psi|^2,
\end{eqnarray*}
and in the similar way
$$
\text{Re}\langle \nabla
_{e_2}\psi,e_1\cdot\nabla_{e_1}\psi\rangle=-\frac 12
e^u\nabla_{e_2}|\psi|^2.
$$
We also compute
\begin{eqnarray*}
&&\text{Re}\langle
\psi,\nabla_{e_1}(\slashiii{D}\psi)\rangle-i\text{Re}\langle
\psi,\nabla_{e_2}(\slashiii{D}\psi)\rangle\\
&=&-\text{Re}\langle
\psi,\nabla_{e_1}(e^u\psi)\rangle+i\text{Re}\langle
\psi,\nabla_{e_2}(e^u\psi)\rangle\\
&=&-2e^u|\psi|^2\partial_zu-e^u\partial_z|\psi|^2
\end{eqnarray*}

\noindent Therefore we get
\begin{equation*}
\partial_{\bar{z}} T_2(z)= e^u\partial_z|\psi|^2-e^u|\psi|^2\partial_zu.
\end{equation*}
\noindent Hence
$$
\partial_{\bar{z}}T(z)=\partial_{
\bar{z}}T_1(z) +\frac 14\partial_{\bar{z}}T_2(z)=-\frac
14\partial_z K_g.
$$

Therefore $\partial_{\bar{z}}T(z)=0$ when $K_g$ is constant and
$T(z)$ is holomorphic.  We finish the proof.

\end{proof}

\begin{rem}
It is well-known that every holomorphic quadratic differential on
$\S^2$ vanishes identically (see \cite{Jo}). Therefore $T(z)=0$ if
$M=\S^2$.
\end{rem}

\begin{rem}By a similar method as in \cite{CJLW}, we can construct the
holomorphic quantity in the following way. Let ($u,\psi$) be a
solution of (\ref{1}) on $M$. Define a tensor
\begin{eqnarray*}
T_{\alpha \beta }&=&2(u_{_\alpha },u_{_\beta })-\delta _{\alpha
\beta }\sum_{r=1}^2(u_r,u_r)-2u_{\alpha \beta }+\delta _{\alpha
\beta }\sum_{r=1}^2u_{rr}+2\text{Re}\langle \psi ,e_{_\alpha
}\cdot
\nabla _{e_{_\beta }}\psi \rangle\\
& & +\delta _{\alpha \beta }e^u\left| \psi \right| ^2
\end{eqnarray*}
where $u_{_\alpha }=\nabla _{e_\alpha }u,$ and \{$e_1,e_2\}$ is a
local orthonormal basis on $M$. Then we can check  as in Proposition \ref{pro2},
\begin{enumerate}
\item  $T_{11}+T_{22}=0,$

\item  $T_{\alpha \beta }=T_{\beta \alpha },$  i.e., the tensor
$T_{\alpha \beta }$ is symmetric.

\item  $\sum_{\alpha =1}^2\nabla _{e_{\alpha }}T_{\alpha \beta
}=-\partial_{\beta}K_g$.
\end{enumerate}
\par

\noindent  Define $T(z)=\frac 14(T_{11}-iT_{12})$. Then $T(z)dz^2$
is the  holomorphic quadratic differential of Proposition
\ref{pro2}.
\end{rem}

\section{Compactness Theorem}

In this section we consider the compactness of solutions of (
\ref{1}) under the condition that
\[I(u,\psi):=\int_M (e^{2u}+ |\psi|^4) dv<C.\]
Since (\ref{1}) is conformally invariant, in general the set of
solutions of (\ref{1}) with a uniformly bounded energy $I(u,\psi)$
is non-compact.

First, we define weak solutions of (\ref{1}). We say that
$(u,\psi)$ is a weak solution of (\ref{1}), if $u\in W^{1,2}(M)$
and $\psi \in W^{1,\frac 43}(\Gamma (\Sigma M))$ satisfy
\begin{eqnarray*}
\int_M \nabla u\nabla\phi dv&=& \int_M (2e^{2u}-e^u
|\psi|^2-K_g)\phi dv\\
\int_M \langle \psi,\slashiii{D}\xi\rangle
dv&=&-\int_Me^u\langle\psi,\xi\rangle dv
\end{eqnarray*}
for any smooth function $\phi$ and any smooth spinor $\xi$.
It is clear that
$(u,\psi)\in
 W^{1,2}(M)\times W^{1,\frac 43}(\Gamma (\Sigma M))$ is a weak solution
 if and only if $(u,\psi)$ is a critical point of $E$ in
  $W^{1,2}(M)\times W^{1,\frac 43}(\Gamma (\Sigma M))$. A weak solution
  is a classical solution by the following

\begin{prop}\label{prop-regu}
Any weak solution $(u, \psi)$  to (\ref{1}) on $M$ with
$I(u,\psi)<\infty $ is smooth.
\end{prop}

To prove the proposition, we first need a basic inequality in
\cite{BM}.

\begin{lm}\label{lmbm}
Assume $\Omega \subset \R^2$ is a bounded domain and let $u$ be a
solution of
\[
\left\{
\begin{array}{rcll}
\ds\vs -\Delta u &=&\ds f(x)\qquad &\hbox{ in }\Omega  \\
u&= &0\qquad \qquad \,\,&\hbox{ on }\partial \Omega
\end{array}
\right.
\]
with $f\in L^1(\Omega ).$ Then for every $\delta \in (0,4\pi )$ we have
\begin{equation}
\int_\Omega \exp \{\frac{\left( 4\pi -\delta \right) \left| u(x)\right| }{%
\left\| f\right\| _1}\}dx\leq \frac{4\pi ^2}\delta (\text{diam}\Omega )^2,
\label{4}
\end{equation}
where $\left\| f\right\| _1=\int_\Omega \left| f(x)\right| dx.$
\end{lm}

Let $B_r=B_r(x)$ be a geodesic ball at a point $x$ on $M$ with
radius $r$. Here $r$ is smaller than the injective radius of $M$.

\begin{lm}\label{lmint}
If $(u,\psi )$ is a weak solution to (\ref{1}) in $B_r$ satisfying  $
\int_{B_r}e^{2u}+|\psi|^4dx<\infty$, then we have
$$
u^+ \in L^{\infty}(B_{\frac r4})\quad \hbox{ and } \quad \left|\psi\right|\in
L^{\infty}(B_{\frac r4}).
$$
\end{lm}

\begin{proof}
First we consider  $u.$  Set
\[
f_1=2e^{2u}-e^u\left| \psi \right| ^2-K_g.
\]
Then we have
\[
-\Delta u=f_1.
\]
We consider the following Dirichlet problem
\begin{equation}
\left\{
\begin{array}{rcl}
-\Delta u_1 &=& f_1,\qquad \text{in  }B_r \\
 u_1&=& 0,
\qquad \text{ on  }\partial B_r.
\end{array}
\right.  \label{57}
\end{equation}
Since $\int_{B_r}e^{2u}dx<\infty$ and $\int_{B_r}\left| \psi
\right| ^4dx<\infty$ we know that $f_1\in L^1(B_r).$ By
applying Lemma
\ref{lmbm} on a smaller domain we have
\begin{equation}
e^{k\left| u_1\right| }\in L^1(B_r)  \label{66}
\end{equation}
for some $k>1$ and in particular $ u_1\in {L^p(B_r)}$ for
some $p> 1.$

Let $u_2=u-u_1$ so that $\Delta u_2=0$ on $B_r.$ The mean value
theorem for harmonic functions implies that
\[
\left\| u_2^{+}\right\| _{L^\infty (B_{\frac r2})}\leq C\left\|
u_2^{+}\right\| _{L^1(B_r)}.
\]
Since $u_2^{+}\leq u^{+}+\left| u_1\right| $ and
$2\int_{B_r}u^{+}\leq \int_{B_r}e^{2u}<\infty,$ we have
$u_2^{+}\in{L^1(B_r)}$ and consequently

\begin{equation}
\left\| u_2^{+}\right\| _{L^\infty (B_{\frac r2})}<\infty.
\label{55}
\end{equation}
Next we write
\[
f_1=2e^{2u_2}e^{2u_1}-e^{u_1}e^{u_2}\left| \psi \right| ^2-K_g.
\]
 From (\ref{55}) and (\ref{66}) we have $f_1\in L^{1+\varepsilon
}(B_{\frac r2})$ for some $\varepsilon
>0.$ Hence standard elliptic estimates imply that

\[
\left\| u^{+}\right\| _{L^\infty (B_{\frac r4})}\leq C\left\|
u^{+}\right\| _{L^1(B_r)}+C\left\| f_1\right\| _{L^{1+\varepsilon
}(B_{\frac r2})}<\infty .
\]

Since $u^+ \in L^{\infty}(B_{\frac r4})$, then the right hand of
equation $\slashiii{D} \psi =-e^u\psi$ is in $L^4(\Gamma (\Sigma
B_{\frac r4}))$. Hence $\psi \in C^0(\Gamma (\Sigma B_{\frac
r4}))$ and especially $|\psi| \in L^{\infty}(B_{\frac r4})$.
\end{proof}

\noindent{\it Proof of Proposition \ref{prop-regu}.}
The standard method, together
with Lemma 4.3,
implies that  $u$ and $\psi$ are smooth.
\qed

\

Next we discuss the compactness of a sequence of smooth solutions to (\ref{1}).
We begin with studying uniformly $L^\infty $ boundedness of solutions for (\ref{1}%
). Assume that ($u_n,\psi _n)$ is a sequence of solutions of
(\ref{1}). Similarly as before we set
$$
f_{1}^n =2e^{2u_n}-e^{u_n}\left| \psi _n\right| ^2-K_g,
$$

\begin{lm}\label{lmuni}
Let $\varepsilon _0<\pi$ be a constant.  For any
 sequence of solutions $(u_n,\psi_n)$ with
\[
\int_{B_r}e^{2u_n}dx<\varepsilon _0,\qquad \int_{B_r}\left| \psi
_n\right| ^4dx<C
\]
for some fixed constant $C>0$ we have that $\left\| u_n^{+}\right\|_{L^{\infty}
(B_{\frac r4})}$ is uniformly bounded.
\end{lm}

\begin{proof}
Similarly as in the proof of lemma \ref{lmint}, it is sufficient to
show that $f_{1}^n$ is uniformly bounded in $L_{loc}^q(B_r)$ for
some $q>1$.

Let $w_n$ be the solution of following problem:
\begin{equation*}
\left\{
\begin{array}{rcll}
-\Delta w_n&=&2e^{2u_n},&\qquad  \text{in} \qquad  B_r(x)\\
w_n&=& 0, &\qquad \text{on} \qquad \partial B_r(x).
\end{array}
\right.
\end{equation*}
It is clear that $w_n\geq 0$ in $B_r(x)$.  Since $\varepsilon_0
<\pi$, we can choose $\delta >0$ such that $4\pi -\delta
>2\varepsilon_0(2+\delta)$. By lemma \ref{lmbm} we get
\begin{equation}\label{aa}
\int_{B_r(x)}e^{(2+\delta)w_n}\leq C
\end{equation}
for some constant $C$.

Next let $z$ be the solution of  the following equation
\begin{equation*}
\left\{
\begin{array}{rcll}
-\Delta z &=&-K_g,& \qquad  \text{in}  \qquad B_r(x)\\
z&=&0, &\qquad \text{on} \qquad \partial B_r(x).
\end{array}
\right.
\end{equation*}

\noindent It is clear that $\Delta (u_n-w_n-z)=e^u|\psi|^2\geq 0$
on $B_r(x)$ and
\begin{eqnarray*}
\int_{B_r(x)}(u_n-w_n-z)^+&\leq & \int_{B_r(x)}(u_n-w_n)^++|z|dx \\
&\leq & \int_{B_r(x)}(u_{n}^{+}+|z|)
\leq \int_{B_r(x)}e^{2u_n}+C_1 \le C,
\end{eqnarray*}
for some constant $C>0$.
Here we have used $w_n\ge 0$.
Therefore, by the mean value theorem for subharmonic function, for
any $y\in B_{\frac r2}(x)$, we have
\begin{eqnarray}\label{ab}
(u_n-w_n-z)(y)&\leq & C\int_{B_r(x)}(u_n-w_n-z)\nonumber\\
&\leq & C\int_{B_r(x)}(u_n-w_n-z)^+\leq C
\end{eqnarray}

\noindent Thus, from (\ref{aa}) and (\ref{ab}), we deduce that
\begin{equation}\label{ac}
\int_{B_{\frac r2}(x)}e^{(2+\delta)u_n}\leq C.
\end{equation}
By the H\"older inequality, for $l=\frac{4+2\delta}{4+\delta}>1$ we
have

$$
\int_{B_{\frac r2}}(e^{u_n}\left| \psi_n \right|^2)^ldx\leq
(\int_{B_{\frac r2
}}e^{(2+\delta)u_n}dx)^{\frac{l}{2+\delta}}(\int_{B_{\frac r2
}}\left|\psi_n\right|^4dx)^ {\frac{2+\delta-l}{2+\delta}}\leq C.
$$
Let $q=\min \{l, 2+\delta\}$. We have established that $f_{1}^n$ is
uniformly bounded in $L^q(B_{\frac r2}(x))$ with $q>1$.
\end{proof}

Since $\left\| u_n^{+}\right\|_{L^{\infty} (B_{\frac r4})}$ is
uniformly bounded, by the standard method and the bootstrapping method
of elliptic equations, we can get  uniform estimates for higher
 derivatives of the functions $u_n$ and $\psi_n$. That is,

\begin{thm}\label{th-re}
Assume that $(u_n, \psi _n)$ is a sequence of solutions for
(\ref{1}) with
\[
\int_{M}e^{2u_n}dv<\varepsilon _0,\text{ and }\int_{M}\left| \psi
_n\right| ^4dv<C
\]
for some positive constant $\varepsilon _0<\pi$ and $C$. Then we
have
\begin{equation}\label{ah}
\| u_n\|_{C^k(B_{\frac 18(x)})}+\| \psi_n\|_{C^k(B_{\frac
18(x)})}\leq C.
\end{equation}
for any geodesic ball $B_{\frac 18(x)}$ of $M$.
\end{thm}

\

>From Theorem \ref{th-re}, we have the following Theorem:
\begin{thm}
Assume that $(u_n,\psi _n)$ is a sequence of solutions for
(\ref{1}) with
\[
\int_{M}e^{2u_n}dv<\varepsilon _0,\text{ and }\int_{M}\left| \psi
_n\right| ^4dv<C
\]
for some positive constant $\varepsilon _0<\pi$ and $C$. Then
$(u_n,\psi _n)$ admits a subsequence converging to $(u,\psi)$
which is a smooth solution of (\ref{1}).
\end{thm}

\section{Blow up behavior}

When the energy $\int_M e^{2u}dv $ is large, then the blow-up
phenomenon may occur as in the case of the Liouville equation. In
this section we will analyze the asymptotic behavior of a sequence
of solutions for (\ref{1}) when the blow-up phenomenon happens.
Assume that ($u_n,\psi _n)$ satisfies
\begin{equation}
\left\{
\begin{array}{rcl}
-\Delta u_n&=&2e^{2u_n}-e^{u_n}\left|\psi _n \right|^2 -K_g,\qquad \\
\slashiii{D}\psi _n&=&-e^{u_n}\psi _n,
\end{array}
\text{on }M\right.  \label{77}
\end{equation}
with
\begin{equation}
\int_{M}e^{2u_n}dv<C,\text{ and }\int_{M}\left| \psi _n\right|
^4dv<C \label{78}
\end{equation}
for some positive constant $C.$

We shall follow  \cite{BM}, where
the authors analyze the behavior of a sequence of solutions for
the  Liouville-type equation on a bounded domain.  Similar results for the Toda system,
which is another  natural generalization of the Liouville equation,
were obtained in \cite{JW}.

\begin{thm}\label{mainthm}
Let $(u_n,\psi _n)$ be a sequence of solutions to $\left( \ref{77}\right) $
 satisfying $\left( \ref{78}\right) .$ Define
\[\begin{array}{rcl}
\Sigma _1&=&\left\{ x\in M,\text{ there is a sequence
}y_n\rightarrow x\text{ such that }u_n(y_n)\rightarrow +\infty
\right\}\\
\Sigma _2&=&\left\{ x\in M,\text{ there is a sequence
}y_n\rightarrow x\text{ such that }\left| \psi _n(y_n)\right|
\rightarrow +\infty \right\}
.\end{array}
\]
Then, we have $\Sigma_2\subset\Sigma_1$. Moreover,
$(u_n,\psi _n)$ admits a subsequence,
denoted still by $(
u_n,\psi _n),$ satisfying that

\begin{enumerate}
\item[a)] $\psi _n$ is bounded in $%
L_{loc}^{\infty} (M\backslash \Sigma _2)$ .

\item[b)]  For $u_n$, one of the following alternatives holds:
\begin{enumerate}
\item[i)]  $u_n$ is bounded in $L^{\infty} (M).$

\item[ii)]  $u_n$ $\rightarrow -\infty $ uniformly on $M$.

\item[iii)]  $\Sigma _1$ is finite, nonempty and either
\begin{equation}
u_n\text{ is bounded in }L_{loc}^\infty (M\backslash \Sigma _1)
\label{aaa}
\end{equation}
or
\begin{equation}
u_n\rightarrow -\infty \text{ uniformly on compact subsets of
}M\backslash \Sigma _1.  \label{ddd}
\end{equation}
\end{enumerate}
\end{enumerate}
\end{thm}

\begin{proof}
First, if $x\in M\backslash \Sigma _1,$ then from the equation
$\slashiii{D}\psi _n=-e^{u_n}\psi _n$ we know $x\in M\backslash
\Sigma _2.$ Therefore we
have $\Sigma _2\subset \Sigma _1$ and $\psi _n$ are bounded in $%
L_{loc}^\infty (M\backslash \Sigma _2).$

Next let $f^{n}_{1}$ be as before. Since $e^{2u_{n}}$ is bounded in $%
L^1(M)$, we may extract a subsequence from $u_n$ (still denoted
$u_{n}$) such that $e^{2u_{n}}$ converges in the sense of measures
on $M$ to some nonnegative bounded measure $\mu$ i.e.
\begin{equation*}
\int_{M}e^{2u_n}\varphi dv \rightarrow \int_{M}\varphi d\mu
\end{equation*}
for every $\varphi \in C(M).$ A point $x\in M$ is called an
$\varepsilon - $regular point with respect to $\mu$ if there is a
function $\varphi \in C(M),$ $\text{supp}\varphi \subset
B_r(x)\subset M,$ $0\leq \varphi \leq 1$ with $\varphi =1$ in a
neighborhood of $x$ such that
\begin{equation*}
\int_{M}\varphi d\mu  <\varepsilon
\end{equation*}

We define
\begin{equation*}
\Omega _1(\varepsilon ) = \{x\in M:x\text{ is not an }\varepsilon -\text{%
regular point with respect to }\mu\}. \\
\end{equation*}
By definition and (\ref{78}) we see that $\Omega _1(\varepsilon )$ is
finite. We divide the proof into three steps.

\

\noindent{\it Step 1}. $\Sigma _1=\Omega _1(\varepsilon _0)$
provided $\varepsilon _0<\pi.$

First we show that $\Omega _1(\varepsilon _0)\subset \Sigma _1.$
Supposing that $x_0\in \Omega _1(\varepsilon _0),$ we claim that
for any $R>0,$ $\lim_{n\rightarrow +\infty }\left\|
u_n^{+}\right\|_{L^\infty (B_R(x_0))}=+\infty .$ We demonstrate
the claim by a contradiction. So we assume that there would be
some $R_0>0$ and a subsequence such that $\left\|
u_n^{+}\right\| _{L^\infty (B_{R_0}(x_0))}$ is bounded. Especially we have $%
\left\| e^{2u_n}\right\| _{L^\infty (B_{R_0}(x_0))}$ $\leq C$ and therefore $%
\int_{B_R(x_0)}e^{2u_n}dx\leq CR^{\delta}$ for all $R<R_0$ and some $%
\delta >0.$ This implies
\[
\int_{M}\varphi d\mu <\varepsilon _0\text{ for some suitable
}\varphi .
\]
Therefore $x_0$ is regular, contradicting $x_0\in \Omega
_1(\varepsilon _0).$ So the claim is proved. Now we choose $R>0$
small enough so that $\overline{B}_R(x_0)$ does not contain any
other point of $\Omega _1(\varepsilon _0).$ Let $x_n\in B_R(x_0)$
be such that
\[
u_n^{+}(x_n)=\max_{\overline{B}_R(x_0)}u_n^{+}\rightarrow +\infty .
\]
We claim that $x_n\rightarrow x_0,$ i.e. $x_0\in \Sigma _1.$ Otherwise there
would be a subsequence
\[
x_{n_k}\rightarrow \overline{x}\neq x_0\text{ and }\overline{x}\notin \Omega
_1(\varepsilon _0)
\]
that is, $\overline{x}$ is a regular point. This is a
contradiction. Therefore we have proved that $\Omega
_1(\varepsilon _0)\subset \Sigma _1.$

Next we show that $\Sigma _1\subset \Omega _1(\varepsilon _0)$ by
using the approach to the Toda system in \cite{JW}. Let $x_0\in \Sigma
_1.$ Assume by
contradiction that $x_0\notin \Omega _1(\varepsilon _0).$ Thus $%
\int_{B_\delta \left( x_0\right) }e^{2u_n}\leq \varepsilon _0$ for
any small constant ${\delta}>0.$ Note that $-\Delta
u_n=2e^{2u_n}-e^{u_n}\left|\psi _n\right|^2-K_g\leq
2e^{2u_n}-K_g.$ Define $w:B_\delta (x_0)\rightarrow R$ by
\begin{equation}
\left\{
\begin{array}{rcl}
-\Delta w &=& \ds\vs  2e^{2u_n}-K_g\,\qquad\text{in }B_\delta (x_0) \\
w &=& u_n\, \qquad\text{on }\partial B_\delta (x_0).
\end{array}
\right.
\end{equation}
The maximum principle implies that $u_n\leq w.$ Since $\Sigma _1$
is finite, we may assume that $u_n$ is uniformly bounded in
$L^\infty (\partial B_\delta (x_0)).$ In view of $\int_{B_\delta
\left( x_0\right) }e^{2u_n}\leq \varepsilon _0<\pi,$ and the
boundedness of the  curvature $R$ of $M$,  as in the proof
of lemma \ref{lmuni} we have $w^{+}\in L^\infty (B_{\frac \delta
2}(x_0)),$ which in turn implies that $u_n^{+}\in L^\infty
(B_{\frac \delta 2}(x_0)).$ Hence we have a contradiction.
Therefore $\Sigma _1\subset \Omega _1(\varepsilon _0).$

So we have $\Sigma _1=\Omega _1(\varepsilon _0).$

\

\noindent{\it Step 2.} $\Sigma _1=\emptyset $ implies (1) and (2) hold.

$\Sigma _1=\emptyset $ means that $u_n^{+}$ is bounded in
$L^\infty
(M).$ Consequently $\psi _n$ is bounded in $L^\infty (M).$ Thus, $
f^{n}_{1}$ is bounded in $L^p(M)$ for any $p>1$.  Applying the
Harnack inequality as in \cite{BM}, we have (1) or (2).

\

\noindent{\it Step 3.} $\Sigma _1\neq \emptyset $ implies (3).

In this case, we know that $u_n^{+}$ is bounded in $L_{loc}^\infty
(M\backslash \Sigma _1)$ and therefore $f^{n}_{1}$ is bounded in $%
L_{loc}^p(M\backslash \Sigma _1)$ for any $p>1.$ Then as in step (2)
we know that either
\[
u_n\text{ is bounded in }L_{loc}^\infty (M\backslash \Sigma _1),
\]
or
\[
u_n\rightarrow -\infty \text{ \qquad on any compact subset of
}M\backslash \Sigma _1
\]
Thus we complete the proof of the Theorem.%
\end{proof}

Actually in Theorem \ref{mainthm} the case (\ref{aaa}) will not
occur if $\Sigma _1 \backslash \Sigma_2 \neq \emptyset$. Next
we will show this.

\begin{thm}\label{mainthm1}
In Theorem \ref{mainthm}, if in addition $\Sigma _1 \backslash
\Sigma_2 \neq \emptyset$, then the first case of $iii)$ does not
happen, i.e. $ u_n\rightarrow -\infty $ uniformly on compact
subsets of $M\backslash \Sigma _1.$ Moreover, setting $\Sigma
_1=\{p_1,p_2,\cdots ,p_l\},$ we have
\[
e^{2u_n}\rightharpoonup \mu =\sum_{i=1}^l\alpha _i\delta _{p_i},
\qquad \text{with} \quad \alpha_i \geq \pi.
\]
\end{thm}

\begin{proof}
We should show that (\ref{aaa}) does not happen when $\Sigma _1
\backslash \Sigma_2 \neq \emptyset $. Fix some point $x_0\in
\Sigma _1 \backslash \Sigma_2$ and choose $\delta
>0$ to be so small  that $x_0$ is the only point of $\Sigma _1 \backslash \Sigma_2$ in
$\overline{B}_{\delta} (x_0)$. Let $f^{n}_{1}$ be as before, i.e.
$$
f^{n}_{1} =2e^{2u_n}-e^{u_n}\left| \psi _n\right| ^2-K_g.
$$
Since $x_0$ is a point of $\Sigma _1 \backslash \Sigma_2$, we can
select $\delta$ to be sufficiently small such that
\begin{eqnarray*}
f^{n}_{1}&=& e^{2u_n}(2- e^{-u_n}|\psi_n|^2-K_g
e^{-2u_n})\\
&=& e^{2u_n}v_n(x),
\end{eqnarray*}
where $v_n(x)=2- e^{-u_n}|\psi_n|^2-K_g e^{-2u_n}$ and
$v_n(x)\rightarrow 2$ in $B_{\delta}(x_0)$. Therefore we can
rewrite the first equation of (\ref{77}) as
\begin{equation*}
\left\{
\begin{array}{lcl}
-\Delta u_n = v_n(x)
e^{2u_n}, &&\qquad\text{in }B_{\delta}(x_0)\\
0 \leq v_n(x) \leq b,&& \qquad\text{in }B_{\delta}(x_0)\\
\int_{B_{\delta}(x_0)}e^{2u_n} dx \leq C
\end{array}
\right.
\end{equation*}
for $ b$ and $C$ positive constants.

Noting that $x_{0}$ is a blow up point for $u_n$, we can apply the
Brezis-Merle result (see \cite{BM}) to conclude that
$$
u_n \rightarrow -\infty,  \qquad \text{for any compact subset}
\qquad  K\subset B_{\delta}(x_0) \backslash \{x_0\}.
$$
Consequently, by the alternative proved in Theorem \ref{mainthm},
we have that (\ref{aaa}) does not happen and only (\ref{ddd})
holds.

Moreover since the case $\left( \ref{ddd}\right) $ is valid, then
$e^{2u_n}\rightarrow 0 $ in $L_{loc}^p(M\backslash \Sigma _1)$ for
any $p\geq 1$. Therefore, if $e^{2u_n}\rightharpoonup \mu $, then
the measure $ \mu $ is supported on $\Sigma _1.$ Hence, setting
$\Sigma _1=\{p_1,p_2,\cdots ,p_l\},$ we have
$e^{2u_n}\rightharpoonup \mu =\sum_{i=1}^l\alpha _i\delta _{p_i}$
with $\alpha_i \geq \pi$.
\end{proof}

\section{Asymptotic behavior of rescaling equations}

It is well known that a  ``bubble'',  an entire solution of
(\ref{1}) with finite energy, will be obtained after a suitable
rescaling at a blow-up point. In the rest of the paper we will
analyze the asymptotic behavior of an entire solution with finite
energy. We will show that an entire solution on $\R^2$ can be
extended to $\S^2$, i.e. the singularity at infinity is removable.

The considered equations are
\begin{equation}\label{se}
\left\{
\begin{array}{rcl}
-\Delta u &=& \ds\vs 2e^{2u}-e^u\left\langle \psi ,\psi
\right\rangle,\qquad x\in \R^2
\\
\slashiii{D}\psi &=&\ds  -e^u\psi, \qquad x\in \R^2.
\end{array}
\right.
\end{equation}

\noindent The energy condition is
\begin{equation}\label{sec}
I(u,\psi)=\int_{\R^2}(e^{2u}+|\psi|^4)dx<\infty.
\end{equation}

\ \

Next we start to deal with the asymptotic behavior of solutions of
(\ref{se}) and (\ref{sec}). First we have

\begin{lm}\label{blm}
Let $(u,\psi)$ be a solution of (\ref{se}) and (\ref{sec}) with
$u\in H^{1,2}_{loc}(\R^2)$ and $\psi\in H^{1,\frac
43}_{loc}(\R^2)$. Then $u^{+}\in L^{\infty}(\R^2)$.
\end{lm}

The proof of Lemma \ref{blm} follows from the idea of \cite{CL2}.
Since $u^{+}\in L^{\infty}(\R^2)$, it follows from the discussion
in the previous section that $(u,\psi)$ is smooth in $\R^2$.

\
\

Denote  $(v,\phi)$ be the Kelvin transformation of $(u,\psi)$,
i.e.
\begin{eqnarray*}
& & v(x)=u(\frac {x}{|x|^2})-2\ln |x|\\
& & \phi (x)=|x|^{-1} \psi (\frac{x}{|x|^2})
\end{eqnarray*}
Then $(v,\phi)$ satisfies
\begin{equation}\label{sek}
\left\{
\begin{array}{rcll}
-\Delta v &=& 2e^{2v}-e^v\left\langle \phi ,\phi
\right\rangle,&\qquad x\in \R^2\backslash\{0\}
\\
\slashiii{D}\phi &=&\ds  -e^v\phi, &\qquad x\in
\R^2\backslash\{0\}.
\end{array}
\right.
\end{equation}
And, by change of variable,
\begin{eqnarray*}
& & \int_{|x|\leq r_0}e^{2v}dx=\int_{|x|\geq \frac
{1}{r_0}}e^{2u}dx\\
& &\int_{|x|\leq r_0}|\phi|^4dx=\int_{|x|\geq \frac
{1}{r_0}}|\psi|^4dx
\end{eqnarray*}
 becomes small if $r_0$ is small. Therefore we obtain that there
is a $r_0$ small enough such that $(v,\phi)$ is a smooth solution
to (\ref{sek}) on $B_{r_0}\backslash\{0\}$ with energy
$\int_{|x|\leq r_0}e^{2v}dx<\varepsilon_0<\pi$ for any
sufficiently small positive number $\varepsilon_0$, and
$\int_{|x|\leq r_0}|\phi|^4dx<C$. Since (\ref{sek}) and
(\ref{sec}) are conformally invariant, in the sequel we may assume
$B_{r_0}$ to be the unit disk $B_1$.

\begin{lm}\label{asy-phi}
There is an $0<\varepsilon_0 < \pi$ if $(v,\phi)$ is a smooth
solution to (\ref{sek}) on $B_1\backslash\{0\}$ with energy
$\int_{|x|\leq 1}e^{2v}dx<\varepsilon_0$, and $\int_{|x|\leq
1}|\phi|^4dx<C$, then for any $x\in B_{\frac {1}{2}}$ we have
\begin{equation}\label{asy-phi1}
|\phi(x)||x|^{\frac 12}+|\nabla\phi(x)||x|^{\frac 32}\leq
C(\int_{B_{2|x|}}|\phi|^4dx)^{\frac 14}.
\end{equation}
Furthermore, if we assume that $e^{2v}=O(\frac
{1}{|x|^{2-\varepsilon}})$, then,  for any $x\in B_{\frac 12}$, we
have
\begin{equation}\label{asy-phi11}
|\phi(x)||x|^{\frac 12}+|\nabla\phi(x)||x|^{\frac 32}\leq
C|x|^{\frac {1}{4C}}(\int_{B_1}|\phi|^4dx)^{\frac 14},
\end{equation}
for some positive constant $C$. Here $\varepsilon$ is any
sufficiently small positive number.
\end{lm}

\begin{proof} We use a similar argument as in \cite{CJLW}
to prove the Lemma. Fix any $x_0\in B_{\frac 12}\backslash\{0\}$,
and define $(\widetilde{v},\widetilde{\phi})$ by
\begin{eqnarray*}
\widetilde{v}(x) &=& v(x_0+|x_0|x)+\log{|x_0|},\\
\widetilde{\phi}(x) &=& |x_0|^{\frac 12}\phi(x_0+|x_0|x).
\end{eqnarray*}
It is clear that $(\widetilde{v},\widetilde{\phi})$ is a smooth
solution to (\ref{se}) on $B_1$ with
$\int_{B_1}e^{2\widetilde{v}}dx<\varepsilon_0$ and
$\int_{B_1}|\widetilde{\phi}|^4dx<C$. Applying Theorem
\ref{th-re}, we have
$$
|\widetilde{\phi}|_{C^1(B_{\frac 12})}\leq
C|\widetilde{\phi}|_{L^4(B_1)}.
$$
Scaling back, we obtain (\ref{asy-phi1}).

\ \

Next recall that the spinor field $\phi(x)$ satisfies
$$
\slashiii{D}\phi=-e^v\phi \qquad \text{in} \quad
B_1\backslash\{0\}.
$$
We choose a cut-off function $\eta_{\varepsilon}\in
C^{\infty}_{0}(B_{2\varepsilon})$ such that $\eta_{\varepsilon}=1$
in $B_{\varepsilon}(0)$ and $|\nabla \eta_{\varepsilon}|<\frac
{C}{\varepsilon}$. Then we have
$$
\slashiii{D}((1-\eta_{\varepsilon})\phi)=-(1-\eta_{\varepsilon})e^v\phi
-d \eta_{\varepsilon}\cdot \phi.
$$
>From the elliptic estimate with boundary (see \cite{CJLW}), we
have
\begin{eqnarray}\label{asy-phi2}
& & ||(1-\eta_{\varepsilon})\phi||_{W^{1,\frac 43}(B_1)} \nonumber\\
&\leq & C ||e^{v}||_{L^2(B_1)}||\phi||_{L^4(B_1)}
+C||\phi||_{W^{1,\frac 43}(\partial B_1)}+C||d
\eta_{\varepsilon}\cdot \phi||_{L^{\frac 43}(B_1)}.
\end{eqnarray}
By (\ref{asy-phi1}) we have
$$
\lim_{\varepsilon \rightarrow 0}\frac
1\varepsilon||\phi||_{L^{\frac 43}(B_{2\varepsilon})}=0.
$$
Now letting $\varepsilon \rightarrow 0$, and in virtue of the
smallness of $\int_{B_1}e^{2v}dx$ and the Sobolev embedding
theorem, we obtain
$$
(\int_{B_1}|\phi|^4dx)^{\frac 14}\leq C((\int_{\partial
B_1}|\nabla \phi|^{\frac 43}ds)^{\frac 34}+(\int_{\partial
B_1}|\phi|^4ds)^{\frac 14}).
$$
By rescaling, we have for any $0\leq r\leq 1$
\begin{eqnarray*}
(\int_{B_{r}}|\phi|^4dx)^{\frac 14}& \leq &  C(r\int_{\partial
B_{r}}|\nabla \phi|^{\frac 43}ds)^{\frac 34}+C(r\int_{\partial
B_{r}}|\phi|^4ds)^{\frac 14}\\
&\leq &C(r\int_{\partial B_{r}}|\nabla \phi|^{\frac 43}ds)^{\frac
14}+C(r\int_{\partial B_{r}}|\phi|^4ds)^{\frac 14}.
\end{eqnarray*}
i.e.
\begin{equation}\label{asy-phi3}
\int_{B_{r}}|\phi|^4dx\leq Cr(\int_{\partial B_{r}}|\nabla
\phi|^{\frac 43}ds+\int_{\partial B_{r}}|\phi|^4ds).
\end{equation}

\noindent Next let $\overline{\phi}:=\frac
{1}{|B_1|}\int_{B_1}\phi dx $. Note that
$$
\slashiii{D}(\phi-\overline{\phi})=-e^v(\phi-\overline{\phi})-e^v\overline{\phi}
\qquad \text{in} \quad B_1\backslash\{0\}.
$$
By an similar argument for obtaining (\ref{asy-phi2}) and using
the Poincare inequality, we have
$$
||\phi-\overline{\phi}||_{W^{1,\frac 43}(B_1)}\leq C(
||e^{v}||_{L^2(B_1)}||\phi-\overline{\phi}||_{W^{1,\frac
43}(B_1)}+||\nabla \phi||_{L^{\frac 43}(\partial
B_1)}+||e^v\overline{\phi}||_{L^{\frac 43}(B_1)}).
$$
Again, in virtue of  the smallness of $\int_{B_1}e^{2v}dx$ we
obtain
\begin{eqnarray*}
(\int_{B_1}|\nabla \phi|^{\frac 43}dx)^{\frac 34}&\leq &
C(\int_{\partial B_1}|\nabla \phi|^{\frac 43}ds)^{\frac
34}+C|\overline{\phi}|(\int_{B_1}e^{2v}dx)^{\frac 12}\\
&\leq &C(\int_{\partial B_1}|\nabla \phi|^{\frac 43}ds)^{\frac
34}+C(\int_{B_1}|\phi|^4dx)^{\frac 14
}(\int_{B_1}e^{2v}dx)^{\frac 12}\\
&\leq &C(\int_{\partial B_1}|\nabla \phi|^{\frac 43}ds)^{\frac
34}+\varepsilon_1(\int_{B_1}|\phi|^4dx)^{\frac 34
}+C(\varepsilon_1) (\int_{B_1}e^{2v}dx)^{\frac 34},
\end{eqnarray*}
where $\varepsilon_1$ is a small constant. Hence, for $0\leq r\leq
1$, we have
\begin{equation}\label{asy-phi4}
\int_{B_{r}}|\nabla \phi|^{\frac 43}dx\leq Cr\int_{\partial
B_{r}}|\nabla \phi|^{\frac
43}ds+\varepsilon_1\int_{B_{r}}|\phi|^4dx+C(\varepsilon_1)
\int_{B_{r}}e^{2v}dx.
\end{equation}
Note that $e^{2v}=O(\frac {1}{|x|^{2-\varepsilon}})$ for some
$\varepsilon >0$. We have
\begin{equation}\label{asy-phi5}
\int_{B_r}e^{2v}dx\leq Cr\int_{\partial B_r}e^{2v}ds
\end{equation}
>From (\ref{asy-phi3}),(\ref{asy-phi4}) and (\ref{asy-phi5}), for
any $0\leq r\leq 1$, we obtain
$$
\int_{B_r}e^{2v}dx+\int_{B_{r}}|\nabla \phi|^{\frac
43}dx+\int_{B_{r}}|\phi|^4dx\leq Cr(\int_{\partial
B_r}e^{2v}ds+\int_{\partial B_{r}}|\nabla \phi|^{\frac
43}ds+\int_{\partial B_{r}}|\phi|^4ds),
$$
for some constant $C>0$. Denote
$F(r):=\int_{B_r}e^{2v}dx+\int_{B_{r}}|\nabla \phi|^{\frac
43}dx+\int_{B_{r}}|\phi|^4dx$. Then we get
$$
F(r)\leq CrF'(r).
$$
Integrating this inequality yields
\begin{equation}\label{eq1}
F(r)\leq F(1)r^{\frac 1C}.
\end{equation}

\noindent From (\ref{eq1}),  we can easily get (\ref{asy-phi11}).
Thus we complete the proof of Lemma.
\end{proof}

>From Lemma \ref{asy-phi} and the Kelvin transformation, we obtain
the asymptotic estimate of the spinor $\psi(x)$
\begin{equation}\label{asy-psi1}
|\psi (x)|\leq C|x|^{-\frac 12-\delta_0}\qquad \text{for}\quad |x|
\quad \text{near}\quad \infty
\end{equation}
for some positive number $\delta_0$ provided that $e^{2v}=O(\frac
{1}{|x|^{2-\varepsilon}})$.

\
\

Now let $\alpha =\int_{\R^2}2e^{2u}-e^u|\psi|^2dx$, and a constant
spinor $\xi_0=\int_{\R^2}e^u\psi dx$. It will  turn out that the
constant spinor $\xi_0$ is well defined. Then we have

\begin{prop}\label{asy}
Let $(u,\psi)$ be a solution of (\ref{se}) and (\ref{sec}). Then
$u$ satisfies
\begin{equation}\label{ayu}
u(x)=-\frac {\alpha}{2\pi} \ln{|x|}+C+O(|x|^{-1}) \qquad
\text{for}\quad |x| \quad \text{near}\quad \infty,
\end{equation}

\begin{equation}\label{aypsi}
\psi (x)=-\frac {1}{2\pi}\frac{x}{|x|^2}\cdot
\xi_0+o(|x|^{-1})\qquad \text{for}\quad |x| \quad \text{near}\quad
\infty,
\end{equation}
where $\cdot$ is the Clifford multiplication, $C\in R$ is some
constant, and $\alpha =4\pi$.
\end{prop}

\begin{proof}
First, we analyze  the asymptotic behavior of $u(x)$. To show
(\ref{ayu}), we follow essentially an argument used in \cite{CL1}.
Set
$$
w_1(x)=-\frac{1}{2\pi}\int_{\R^2}(\ln{|x-y|}-\ln{(|y|+1)})e^{2u}dy,
$$
$$
w_2(x)=-\frac{1}{2\pi}\int_{\R^2}(\ln{|x-y|}-\ln{(|y|+1)})e^u|\psi|^2dy.
$$
Then, it is easy to check that
$$
\frac {w_1(x)}{\ln{|x|}}\rightarrow -\frac
{1}{2\pi}\int_{\R^2}e^{2u}dx, \qquad \text{as} \quad
|x|\rightarrow +\infty, \quad \text{uniformly},
$$
$$
\frac {w_2(x)}{\ln{|x|}}\rightarrow -\frac
{1}{2\pi}\int_{\R^2}e^{u}|\psi|^2dx, \qquad \text{as} \quad
|x|\rightarrow +\infty, \quad \text{uniformly},
$$

Moreover, $-\triangle w_1(x)=e^{2u}$ and $-\triangle
w_2(x)=e^u|\psi|^2$ on $\R^2$. Therefore, if we define
$v=u(x)-2w_1(x)+w_2(x)$, we have $\triangle v(x)=0$ on $\R^2$.
Since $u^{+}\in L^{\infty}(\R^2)$ by Lemma \ref{blm}, we get that
$$
v(x)\leq C_1+C_2\ln{|x|},
$$
for $|x|$ sufficiently large, with $C_1$,$C_2$ positive constants.
Therefore, by Liouville's theorem on harmonic functions, $v(x)$
has to be constant and hence we get
$$
\frac {u(x)}{\ln{|x|}}\rightarrow -\frac {\alpha}{2\pi}\qquad
\text{as} \quad |x|\rightarrow +\infty, \quad \text{uniformly}.
$$
Since $\int_{\R^2}e^{2u}dx<+\infty$, the above result implies
$$
\alpha \geq 2\pi.
$$

\
\

Next we show that $\alpha >2\pi$. Assume by contradiction that
$\alpha =2\pi$. Let $(v,\phi)$ be the Kelvin transformation of
$(u,\psi)$. We know $(v,\phi)$ satisfy (\ref{sek}) in
$B_1\backslash\{0\}$. Denote $f(x):=2e^{2v}- e^v|\phi|^2$. Then we
have
$$
-\triangle v=f(x) \qquad \text{in}\quad B_1\backslash\{0\}.
$$
>From the asymptotic estimate (\ref{asy-phi1}), we know that
$f(x)>0$ in a small punctured disk $B_1\backslash\{0\}$. Set
$$
h(x)=-\frac{1}{2\pi}\int_{B_1}\log{|x-y|f(y)dy}
$$
and $g(x)=v(x)-h(x)$. It is clear that $\triangle h=-f$ and
$\triangle g=0$.

\noindent On the other hand, we can check that
$$
\lim_{|x|\rightarrow 0}\frac {v}{-\log{|x|}}=0
$$
which implies
$$
\lim_{|x|\rightarrow 0}\frac
{g(x)}{-\log{|x|}}=\lim_{|x|\rightarrow 0}\frac
{v(x)-h(x)}{-\log{|x|}}=\lim_{|x|\rightarrow 0}\frac {u(\frac
{x}{|x|^2})-2\log{|x|}}{-\log{|x|}}=1.
$$
Since $g(x)$ is harmonic in $B_1\backslash\{0\}$, we have
$g(x)=-\log{|x|}+g_0(x)$ with a smooth harmonic function $g_0$ in
$B_1$. By definition, we have $h(x)>0$. Thus, we have
$$
\int_{B_1}e^{2v}dx=\int_{B_1}e^{2g+h}dx\geq \int_{B_1}\frac
{1}{|x|^2}e^{g_0}dx=+\infty,
$$
which is a contradiction with $\int_{\R^2}e^{2v}dx<\infty$. Hence
we have shown that $\alpha > 2\pi$.

\
\

>From $\alpha >2\pi$, we can improve the estimate for $e^{2u}$ to
\begin{equation}\label{ayu3}
e^{2u}\leq C|x|^{-2-\varepsilon}\qquad \text{for}\quad |x| \quad
\text{near}\quad \infty.
\end{equation}
>From (\ref{ayu3}), and by using  potential analysis we also get
$$
-\frac{\alpha}{2\pi}\ln{|x|}-C\leq u(x)\leq
-\frac{\alpha}{2\pi}\ln{|x|}+C
$$
for some constant $C>0$, see \cite{CL2}.

\noindent Then by using (\ref{ayu3}) and (\ref{asy-psi1}) and
following  the derivation of
gradient estimates in \cite{CK}, we get
\begin{equation*}
|\langle x, \nabla u \rangle+\frac{\alpha}{2\pi}|\leq
C|x|^{-\varepsilon} \qquad \text{for} \quad |x|\quad
\text{near}\quad \infty,
\end{equation*}
consequently we have
\begin{equation}\label{ayu4}
|u_r+\frac{\alpha}{2\pi r}|\leq C|x|^{-1-\varepsilon} \qquad
\text{for} \quad |x|\quad \text{near}\quad \infty.
\end{equation}
In the similar way, we can also get
\begin{equation}\label{ayu5}
|u_{\theta}|\leq C|x|^{-\varepsilon} \qquad \text{for} \quad
|x|\quad \text{near}\quad \infty.
\end{equation}
Here $(r,\theta)$ is the polar coordinate system on $\R^2$ and $C,\
\varepsilon$ are  positive constants. From (\ref{ayu5}) and
(\ref{ayu4}), we can obtain (\ref{ayu}). The idea of proving
(\ref{ayu}) can also be seen in \cite{WZ}.

\
\

Next, we show that $\alpha=4\pi$. Set
$$
T(z)=(\partial_z u)^2-\partial^{2}_{z}u+\frac 14\langle \psi,
dz\cdot\partial_{\bar {z}}\psi\rangle+\frac 14\langle
d\bar{z}\cdot\partial_z\psi,\psi\rangle,
$$
where $\cdot$ is Clifford multiplication. From Proposition
\ref{pro2}, we know that $T(z)$ is a holomorphic function. Using
(\ref{asy-psi1}) and  \ref{ayu}), we have the following expansion
of $T(z)$ near infinity
\begin{eqnarray*}
& & \frac 14 (\frac{\alpha}{2\pi})^2\frac{1}{z^2}-\frac 12
\frac{\alpha}{2\pi}\frac{1}{z^2}+o(\frac
{1}{z^2})+\cdots\\
&=& \frac {1}{2z^2}(\frac 12
(\frac{\alpha}{2\pi})^2-\frac{\alpha}{2\pi})+o(\frac
{1}{z^2})+\cdots
\end{eqnarray*}
Hence, $T(z)$ is a constant and $\frac 12
(\frac{\alpha}{2\pi})^2-\frac{\alpha}{2\pi}=0$, i.e.
$\alpha=4\pi$.

\
\

>From $\alpha =4\pi$, we can improve the estimate for $e^{2u}$ to
\begin{equation}\label{ayu6}
e^{2u}\leq C|x|^{-4} \qquad \text{for} \quad |x|\quad
\text{near}\quad \infty.
\end{equation}
This implies that the constant spinor $\xi_0$ is well defined.

\
\

Finally, we analyze the asymptotic behavior of the spinor $\psi(x)$.
We set
$$
\xi(x)=-\frac{1}{2\pi}\int_{\R^2} \frac {x-y}{|x-y|^2}\cdot
e^u\psi dy,
$$
where $\cdot$ is Clifford multiplication. Since the Green function
of the Dirac operator in $\R^2$ is
$$
G(x,y)=\frac {1}{2\pi}\frac {x-y}{|x-y|^2}\cdot,
$$
for any $x,y\in \R^2$ and $x\neq y$, see \cite{AHM}, we have
$\slashiii{D} \xi=-e^u\psi$.

We compute
\begin{eqnarray}\label{aypsi1}
 |x\cdot \xi(x)-\frac {1}{2\pi}\xi_0| \nonumber
&=& \frac{1}{2\pi}|\int_{\R^2}(\frac{x\cdot(x-y)}{|x-y|^2}+1)\cdot
e^u\psi(y)dy| \nonumber\\
&=& \frac{1}{2\pi}|\int_{\R^2}
(\frac{x\cdot(x-y)}{|x-y|^2}-\frac{(x-y)\cdot(x-y)}{|x-y|^2})\cdot
e^u\psi(y)dy| \nonumber\\
&=& \frac{1}{2\pi}|\int_{\R^2}\frac{(x-y)\cdot y}{|x-y|^2}\cdot
e^u\psi(y)dy| \nonumber\\
&\leq & \frac{1}{2\pi} \int_{\R^2}\frac{|y|}{|x-y|}e^u|\psi|dy.
\end{eqnarray}

>From (\ref{ayu6}), we also have
\begin{equation}
|\psi|e^u\leq C|x|^{-2-\varepsilon} \qquad \text{for} \quad
|x|\quad \text{near}\quad \infty,
\end{equation}
for some positive constants $C$ and $\varepsilon$. Then following
the derivation of gradient estimates in \cite{CK}, we get

\begin{equation}\label{aypsi3}
|x\cdot \xi(x)-\frac{1}{2\pi}\xi_0|\leq C|x|^{-\varepsilon} \qquad
\text{for} \quad |x|\quad \text{near}\quad \infty.
\end{equation}

Set $\eta (x)=\psi(x)-\xi(x)$. Since $\slashiii{D}\psi=-e^u\psi$,
we have $\slashiii{D}\eta(x)=0$. By (\ref{asy-psi1}) and
(\ref{aypsi3}) we have $|\eta (x)|\leq C|x|^{-1-\delta_0}$, which
implies $\eta(x)=0$, i.e. $\psi(x)=\xi(x)$. Furthermore,
\begin{eqnarray*}
|\psi(x)+\frac{1}{2\pi}\frac{x}{|x|^2}\cdot \xi_0|
&=& |\frac{x}{|x|^2}\cdot(x\cdot \psi(x)-\frac{1}{2\pi}\xi_0)|\\
&\leq & \frac{1}{|x|}|x\cdot\psi(x)-\frac{1}{2\pi}\xi_0|\\
&\leq & C|x|^{-1-\varepsilon},
\end{eqnarray*}
for $|x|$ near $\infty$. This proves (\ref{aypsi}).
\end{proof}

\

Since the equation (\ref{se}) is conformally invariant, the
solutions $u$ and $\psi$ of (\ref{se}) can be viewed as a function
and a spinor on $ \S^2\backslash{\{northpole\}}$ with finite
energy. In the following Theorem, we shall prove that such a
singularity can be removed as in many conformal problems. Hence,
at the end we obtain that the solutions are actually defined on
$\S^2$.
\begin{thm}
Let $(u,\psi)$ be a smooth solution of (\ref{se}) and (\ref{sec}).
Then $(u,\psi)$ extends to a smooth solution on $\S^2$.
\end{thm}

\begin{proof}
Let $(v,\phi)$ be the Kelvin transformation of $(u,\psi)$. Then
$(v,\phi)$ satisfies (\ref{sek}) on $\R^2\backslash \{0\}$. To
prove the Theorem, it is sufficient to show that $(v,\phi)$ is smooth
on $\R^2$. Applying Proposition \ref{asy}, we have
\begin{equation}
v(x)=(\frac {\alpha}{2\pi}-2) \ln{|x|}+O(1) \qquad \text{for}\quad
|x| \quad \text{near}\quad 0.
\end{equation}
Since $\alpha =4\pi$, we get that $v$ is bounded near $0$. By
recalling that $\phi$ is also bounded near $0$,  elliptic
theory implies that $(v,\phi)$ is smooth.
\end{proof}

\section*{\bf{Acknowledgements}} The work was carried out when the third
author was visiting the Max Planck Institute for Mathematics in
the Sciences. She would like to thank the institute for the
hospitality and the good working conditions.

\end{document}